\documentclass[preprint, 10pt]{elsarticle}
\usepackage{graphicx}

\makeatletter
\newcommand*\bigcdot{\mathpalette\bigcdot@{.5}}
\newcommand*\bigcdot@[2]{\mathbin{\vcenter{\hbox{\scalebox{#2}{$\m@th#1\bullet$}}}}}
\makeatother
\usepackage[margin=2.5cm]{geometry}
\usepackage{setspace}
\usepackage[percent]{overpic}
\usepackage{lmodern}
\usepackage{amsmath,amssymb}
\usepackage{amsfonts,amsthm}
\usepackage{lineno}
\usepackage{graphicx}
\usepackage{multirow}
\usepackage{bm, bbm}
\usepackage{booktabs}
\usepackage{fancyvrb}
\usepackage{algorithm}
\usepackage{algorithmic}

\usepackage{physics}
\usepackage{xcolor}
\usepackage{pxfonts}
\usepackage[footnotesize,bf]{caption}
\usepackage{subcaption}
\usepackage{mathtools}
\usepackage{hhline}
\usepackage[T1]{fontenc}
\usepackage[colorlinks,linkcolor=red,anchorcolor=blue,citecolor=green]{hyperref}
\usepackage{placeins}
\usepackage{makecell}
\usepackage[capitalize]{cleveref}
\usepackage{comment}

\newcommand{\R}{\mathbb{R}}

\newtheorem{remark}{Remark}[section]

\graphicspath{{figures/}}
\usepackage{subcaption}

\journal{}

\begin{document}
\begin{frontmatter}
	\title{Physics-guided correction for operator learning \\ under model misspecification}
    \author[SHNU]{Lei Ma}
	\author[Imperial]{Nicolas Boull\'e}
	\author[TJ]{Yu-Sen Yang}
	\author[SJTU]{Hao Wu} 
	\author[SHNU]{Ling Guo\corref{cor}}
	\cortext[cor]{Corresponding author.} 
	\ead{lguo@shnu.edu.cn}  
	
	\vspace{0.2cm}
	\address[SHNU]{Department of Mathematics, Shanghai Normal University, Shanghai, China}
	
	\vspace{0.2cm}
	\address[Imperial]{Department of Mathematics, Imperial College London, London, UK}
	
	\vspace{0.2cm}
	\address[TJ]{College of Civil Engineering, Tongji University, Shanghai, China}
	
	\vspace{0.2cm}
	\address[SJTU]{School of Mathematical Sciences, Institute of Natural Sciences, and MOE-LSC, Shanghai Jiao Tong University, Shanghai, China}

\begin{abstract}
	Physics-informed operator learning provides an efficient framework for approximating solution operators of partial differential equations by combining observational data with governing physical laws. However, most existing methods implicitly assume that the prescribed governing equation is accurate. This assumption may fail in practical applications, where model simplifications, missing physical effects, parameter drift, or incomplete constitutive relations can lead to model misspecification. In this work, we propose a physics-guided operator correction framework for learning solution operators under misspecified governing equations. At the operator level, the target mapping is decomposed into a prior operator induced by an approximate physical model and a learnable correction operator that accounts for the remaining discrepancy. We realize the operator model correction using a serial DeepONet architecture, where the first DeepONet provides a solution prediction and the second DeepONet learns an additive correction conditioned on both the input function and the solution prediction. The learned correction is incorporated into the physics residual and trained together with data-consistency constraints, allowing the model to retain useful physical structure while adapting to inaccurate governing equations. Numerical experiments on diffusion-reaction, Burgers, cavity flow, and hyperelastic problems show that the proposed method substantially reduces errors induced by misspecified physics. Additional tests under sparse and noisy observations further demonstrate the robustness of the framework and its ability to provide informative uncertainty estimates through deep ensembles.
\end{abstract}


\end{frontmatter}
\section{Introduction}\label{sec:introduction}
Partial differential equations (PDEs) model the spatiotemporal evolution of physical systems arising in fields like fluid mechanics \cite{batchelor2000introduction}, heat transfer \cite{incropera1996fundamentals}, mass transfer \cite{bird2002transport}, solid mechanics \cite{goodier2016elasticity}, and multi-physics coupling \cite{keyes2013multiphysics}. In recent years, scientific machine learning \cite{raissi2019physics, karniadakis2021physics} has driven a paradigm shift from numerical PDE solvers to operator learning \cite{kovachki2023neural, kovachki2024operator,boulle2024mathematical}, with the aim of developing data-driven reduced order models associated with unknown PDEs with small inference time. Representative approaches include Deep Operator Network (DeepONet) \cite{lu2021learning}, Fourier Neural Operator (FNO) \cite{li2020fourier} along with numerous variants that improve architectural structure, generalization across resolutions and discretizations, and computational efficiency \cite{haghighat2024deeponet, qiu2024derivative, tran2021factorized, li2020multipole}. More importantly, Transformer models \cite{vaswani2017attention} and foundation models \cite{bommasani2021opportunities} have emerged as promising architectures for operator learning, offering enhanced capabilities for modeling long-range dependencies and capturing complex relationships in high-dimensional function spaces. Operator learning methods based on Transformer and foundation model architectures are discussed in \cite{boulle2024mathematical, shih2025transformers, hao2023gnot, sun2025towards}.

Data-driven operator learning typically usually requires a large number of high-fidelity input-output pairs, which could be expensive or difficult to acquire in practical applications. To reduce the need for large training datasets, physics-informed learning has been incorporated into neural training objectives. A canonical example is Physics-Informed Neural Networks (PINNs), which enforce PDE residuals and initial or boundary conditions in the loss function, thereby promoting physical consistency under limited observations \cite{raissi2019physics}. This idea has been combined with operator learning, giving rise to Physics-Informed DeepONets \cite{wang2021learning} and Physics-Informed Neural Operators (PINO) \cite{li2024physics}. These methods learn the solution operators associated with known PDEs while imposing PDE constraints at collocation points by a residual loss, computed using automatic differentiation. Despite their efficiency, physics-informed operator learning techniques rely on the assumption that the governing equations used to construct the residual loss accurately represent the true underlying PDE. In practice, model simplifications, unknown or incomplete physics, parameter drift, and unmodeled source terms frequently lead to model misspecification \cite{acquesta2022model, kaszas2019tipping}. In this case, the resulting residual constraint becomes structurally inconsistent with observations, which induces two adverse effects: (i) the optimization is biased toward solutions that satisfy the wrong physics, yielding systematic errors in the learned solution operator; and (ii) the physics loss may conflict with data-consistency terms, producing gradient interference that deteriorates training stability, slows convergence, or traps optimization in poor local minima.

In this work, we propose a physics-guided operator correction framework for learning reliable solution operators under model misspecification. Specifically, at the level of the differential operator, the true mapping is represented as the sum of a prior operator associated with an approximate governing equation and a learnable correction term accounting for the operator discrepancy (see \cref{fig:correction}). This yields a unified and architecture-independent formulation for correcting misspecified physics in operator learning. To instantiate this framework, we adopt a serial DeepONet realization: the first network provides a baseline prediction of the solution function, and the second network, conditioned on both the input function and the solution prediction, produces an additive correction term. This learned correction is embedded into the residual of the prior model and optimized jointly with the available data-consistency constraints. The resulting method retains the benefits of physics-informed operator learning while alleviating the systematic bias and training instability induced by incorrect physical constraints. More broadly, the proposed framework is not restricted to the correction of misspecified PDE models. Whenever a prior operator prediction is available, the same idea can be used to learn a corrective operator for the remaining discrepancy. In particular, the prior component may arise from an approximate numerical solver, a pretrained neural operator, or a prior physics-based operator model. In this sense, the setting of model misspecification considered in this work serves as a natural and representative instance of a more general operator correction paradigm.

\paragraph{Related works}
The modeling and calibration communities have long studied model discrepancy. Classical Bayesian calibration explicitly represents the systematic gap between a computational model and reality via a discrepancy function, jointly inferred with uncertain parameters \cite{kennedy2001bayesian}, and subsequent work has examined prior modeling and identifiability issues \cite{ling2014selection}. In machine learning and dynamical systems, related directions include missing physics learning and discrepancy modeling, where auxiliary models represent unmodeled terms or systematic residuals to improve predictive accuracy and to disentangle deterministic structure from stochastic effects \cite{ebers2024discrepancy}. Related works have exploited prior model knowledge to derive sample complexity estimates and improve training data distribution in operator learning \cite{boulle2024operator,boulle2023elliptic,boulle2022learning,boulle2023learning}. Another relevant line is data-driven equation discovery (e.g., PDE-FIND and SINDy-type extensions), which seeks sparse representations of governing equations from observations \cite{rudy2017data, champion2019data}. Within the PINN framework, related works include symbolic operator discovery under sparse data \cite{podina2022pinn} as well as methods for correcting model misspecification in the governing physics constraints \cite{zou2024correcting}. 

\paragraph{Main contributions} Our contributions are summarized below:
\begin{enumerate}
	\item We introduce a unified framework for operator correction under model misspecification, which decomposes governing equation into a prior operator response and a learned correction at the operator level.
	
	\item As a concrete realization of this framework, we design a serial DeepONet architecture in which the correction operator depends on both the input function and the predicted solution, allowing for state-dependent discrepancies.
	
	\item We formulate a physics-guided training strategy by embedding the learned correction into the residual of the prior model, thereby reducing systematic bias caused by misspecified physics and enhancing predictive accuracy.

	\item We validate the proposed framework on several benchmark problems involving diverse model misspecifications, and further investigate its robustness and uncertainty quantification performance under sparse and noisy observations.
\end{enumerate}

The rest of the paper is organized as follows. \cref{sec:method} develops the proposed operator correction framework under model misspecification and its DeepONet implementation, including the network architecture, training objective, and training algorithm, which are evaluated in \cref{sec:example} through numerical experiments. Finally, \cref{sec:conclusion} concludes the paper with a discussion of future research directions.

\section{Methodology}\label{sec:method}

We consider a partial differential equation described by $\mathcal{N}[u,v]=0$, where \(u\in\mathcal{U}\) denotes the solution and \(v\in\mathcal{V}\) denotes an input function or parameter that specifies the PDE instance under consideration. In this setting, \(v_i\) may represent, for example, a source term, coefficient field, boundary condition, initial condition, material parameter, or other problem-defining input. The aim of operator learning is to approximate the solution operator $\mathcal{G}:\mathcal{V}\to\mathcal{U}$, from pairs \(\{(v_i,u_i)\}_{i=1}^n\) satisfying $\mathcal{N}[u_i,v_i]=0$ for $1\leq i\leq n$. Here, the functions can either be defined on a bounded spatial domain \(\Omega\subset \mathbb{R}^d\), with suitable boundary conditions enforced in the function space or encoded in the input data in the case of boundary value problems, or on a spatio-temporal domain \(\Omega\times [0,T]\) with \(T>0\) for time-dependent PDEs.

In the standard framework of physics-informed neural operators, one may assume that the prescribed PDEs provide an accurate description of the underlying physical system, and that the governing operator $\mathcal{N}$ is exactly known. Under this assumption, the physical model serves as a highly reliable prior and can be directly incorporated into the learning process. Consequently, both data-driven and physics-informed strategies can be employed to train neural operators for approximating the solution operator, typically by minimizing a combination of data misfit and PDE residual losses.

In this work, we consider a more general and practically relevant setting in which the prior PDE model may be structurally imperfect, containing modeling errors, unresolved effects, or missing terms. Such discrepancies are particularly common in complex systems, where deriving an exact governing equation is often challenging due to unresolved effects, unknown interactions, or unavoidable modeling approximations. Enforcing an imperfect PDE model as an exact constraint may introduce systematic bias and restrict the learned operator to an inaccurate solution manifold. Instead of assuming that the governing operator is exactly known, we assume that only a prior operator $\mathcal{N}_0$ is available, which provides an approximate description of the underlying physical system. Therefore, in addition to employing neural operators to approximate the solution operator, we further learn a correction operator to characterize and compensate for the discrepancy $\mathcal N-\mathcal N_0$ between the prior model and the underlying governing operator. By learning this correction operator from observational data, the proposed framework enables adaptive refinement of imperfect physical models while preserving the useful prior knowledge encoded in $\mathcal{N}_0$. In this way, physical information is incorporated as a trustworthy but imperfect prior rather than imposed as an exact constraint, allowing the learned operator to capture latent physical effects and compensate for model deficiencies.

\begin{remark}
	The formulation \(\mathcal{N}[u]=v\) can be recovered as a special case of the above setting, since it can be written equivalently as $\mathcal{N}[u]-v=0$. If \(\mathcal{N}\) is invertible, the induced solution operator satisfies \(\mathcal{G}=\mathcal{N}^{\dagger}\). Thus, the formulation \(\mathcal{N}[u,v]=0\) includes the standard mapping from source terms to solutions as a special case,  while also allowing \(v\) to represent more general problem-defining quantities, such as coefficients, boundary or initial data, and physical parameters.
\end{remark}

\subsection{Operator correction under model misspecification}\label{sec:op-cor}

\begin{figure}[!htbp]
	\centering	
	{\includegraphics[height=0.22\textheight, width=1.0\textwidth]{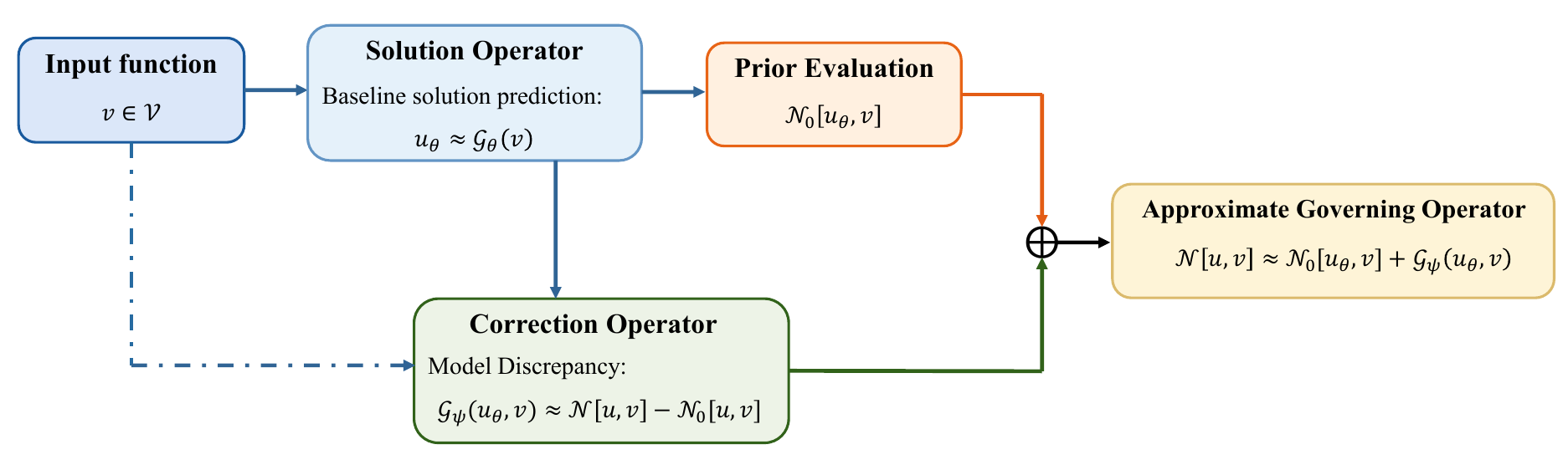}} \\
	\caption{\textbf{Schematic illustration of the unified operator correction framework under model misspecification.} The solution operator produces a baseline prediction $u_{\theta}$, and the correction operator learns the discrepancy conditioned on both the input $v$ and the solution prediction $u_{\theta}$. Their combination constructs an approximation to the true governing operator \(\mathcal{N}[u, v]\) by evaluating the prior operator \(\mathcal{N}_0\) at the learned solution \(u_{\theta}\), and augmenting it with the learned correction \(\mathcal{G}_{\psi}(u_\theta, v)\).}\label{fig:correction}
\end{figure}

Instead of introducing an additive correction directly in the solution space, we formulate the correction at the level of the governing operator. Specifically, we consider the decomposition
\begin{equation*}
	\mathcal{N}\left[u, v\right] = \underbrace{\mathcal{N}_0\left[u, v\right]}_{\text{prior}} + \Bigl(\underbrace{\mathcal{N}\left[u, v\right] - \mathcal{N}_0\left[u, v\right]}_{\text{discrepancy}}\Bigr),\quad u\in\mathcal{U}.
\end{equation*}
This decomposition separates the governing equation into a prior component induced by the approximate physics and a correction term accounting for the effect of model misspecification. We begin by approximating the solution operator by a neural operator $\mathcal{G}_{\theta}:\mathcal{V}\to\mathcal{U}$ such that $u_{\theta} = \mathcal{G}_{\theta}(v)$, where $\theta\in\R^{p_1}$ denotes the trainable parameters of the neural operator. Although the differential operators are written as functions of the solution state, the state $u_{\theta}$ is determined by the problem input $v$. Hence, the operator-level discrepancy evaluated at the prior prediction can depend on both the input $v$ and the predicted solution $u_{\theta}$. We therefore introduce a correction operator $\mathcal{G}_{\psi}$ such that
\begin{equation*}
	\mathcal{G}_{\psi}(u_{\theta}, v) \approx \mathcal{N}\left[u_{\theta}, v\right] - \mathcal{N}_0\left[u_{\theta}, v\right] = \mathcal{N}\left[\mathcal{G}_{\theta}(v), v\right] - \mathcal{N}_0\left[\mathcal{G}_{\theta}(v), v\right],
\end{equation*}
where $\psi\in\R^{p_2}$ denotes the trainable parameters of the correction operator. The corrected operator approximation is then defined by
\begin{equation*}
	\mathcal{N}\left[u, v\right] \approx \mathcal{N}_0\left[u_{\theta}, v\right] + \mathcal{G}_{\psi}(u_{\theta}, v),
\end{equation*}
where \(\mathcal{N}_0\left[u_{\theta}, v\right]\) denotes the prior operator evaluated at the predicted solution, with the differential terms computed by automatic differentiation. This formulation offers two advantages. First, it does not depend on a particular architectural choice: both $\mathcal{G}_{\theta}$ and $\mathcal{G}_{\psi}$ can be represented by any suitable neural operator architectures, such as DeepONet, FNO, etc. Second, it provides a natural justification for allowing the correction to depend on both the input function and the predicted solution. \cref{fig:correction} illustrates the proposed unified operator correction framework. We show how the input \(v\in\mathcal{V}\) is first mapped to a solution prediction \(u_\theta\in\mathcal{U}\), which is then evaluated by the approximate governing operator \(\mathcal{N}_0\) to obtain the prior prediction \(\mathcal{N}_0[u_{\theta}, v]\), and combined with the learned correction \(\mathcal{G}_\psi(u_\theta, v)\) to approximate the true governing operator \(\mathcal{N}[u, v]\).

\subsection{DeepONet realization of the operator correction framework}

In this section, we employ DeepONet as the underlying operator learning model to realize the unified operator correction framework proposed in \cref{sec:op-cor}, which can be referred to as a physics-guided operator correction framework. Specifically, both the solution operator and the correction operator are represented by DeepONet, and the two operators are connected in a serial manner. The first network provides a solution baseline prediction, while the second network further compensates for the discrepancy between the prior model and the true system on the basis of this prediction.

For the solution operator $\mathcal{G}_{\theta}$, which serves as an approximation to the solution function $u$, the corresponding DeepONet representation is given by
\begin{equation*}
	u_{\theta}(y)=\mathcal{G}_{\theta}(v)(y) = \sum_{k=1}^{p_1}\alpha_k\bigl(v(x_1),v(x_2),\dots,v(x_m)\bigr)\beta_k(y),
\end{equation*}
where $\{x_i\}_{i=1}^{m}$ denotes the sensor locations used to sample the input function $v$, and $y=[t,x]^\top$ represents the spatio-temporal coordinate at which the output function is evaluated. The functions $\alpha_k(\cdot)$ and $\beta_k(\cdot)$ correspond to the outputs of the branch and trunk networks, respectively. To realize the correction operator $\mathcal{G}_{\psi}$, we further sample the solution prediction $u_{\theta}$ at a set of sensor locations $\{y_c^{(j)}\}_{j=1}^{n}$. These sampled values, together with the observed input function $v$, are then fed to the correction operator $\mathcal{G}_{\psi}$. The correction operator is then represented by
\begin{equation*}
	\mathcal{G}_{\psi}(u_{\theta}, v)(y) = \sum_{k=1}^{p_2} \gamma_k\bigl(v(x_1),\dots,v(x_m), u_{\theta}(y_c^{(1)}),\dots,u_{\theta}(y_c^{(n)})\bigr)\zeta_k(y),
\end{equation*}
where $\gamma_k(\cdot)$ and $\zeta_k(\cdot)$ correspond to the outputs of the branch and trunk networks of the correction operator, respectively. Under this construction, the final corrected prediction is written as $\mathcal{N}[u](y) \approx \mathcal{N}_0[u_{\theta}(y)]+\mathcal{G}_{\psi}(u_{\theta}, v)(y)$. That is, the first network provides a predicted solution, while the second network yields a correction term, and their sum constitutes the final prediction. In this way, the abstract notion of operator correction introduced in the previous section is realized concretely within the DeepONet framework.

\begin{figure}[t]
	\centering	
	{\includegraphics[height=0.45\textheight, width=1.0\textwidth]{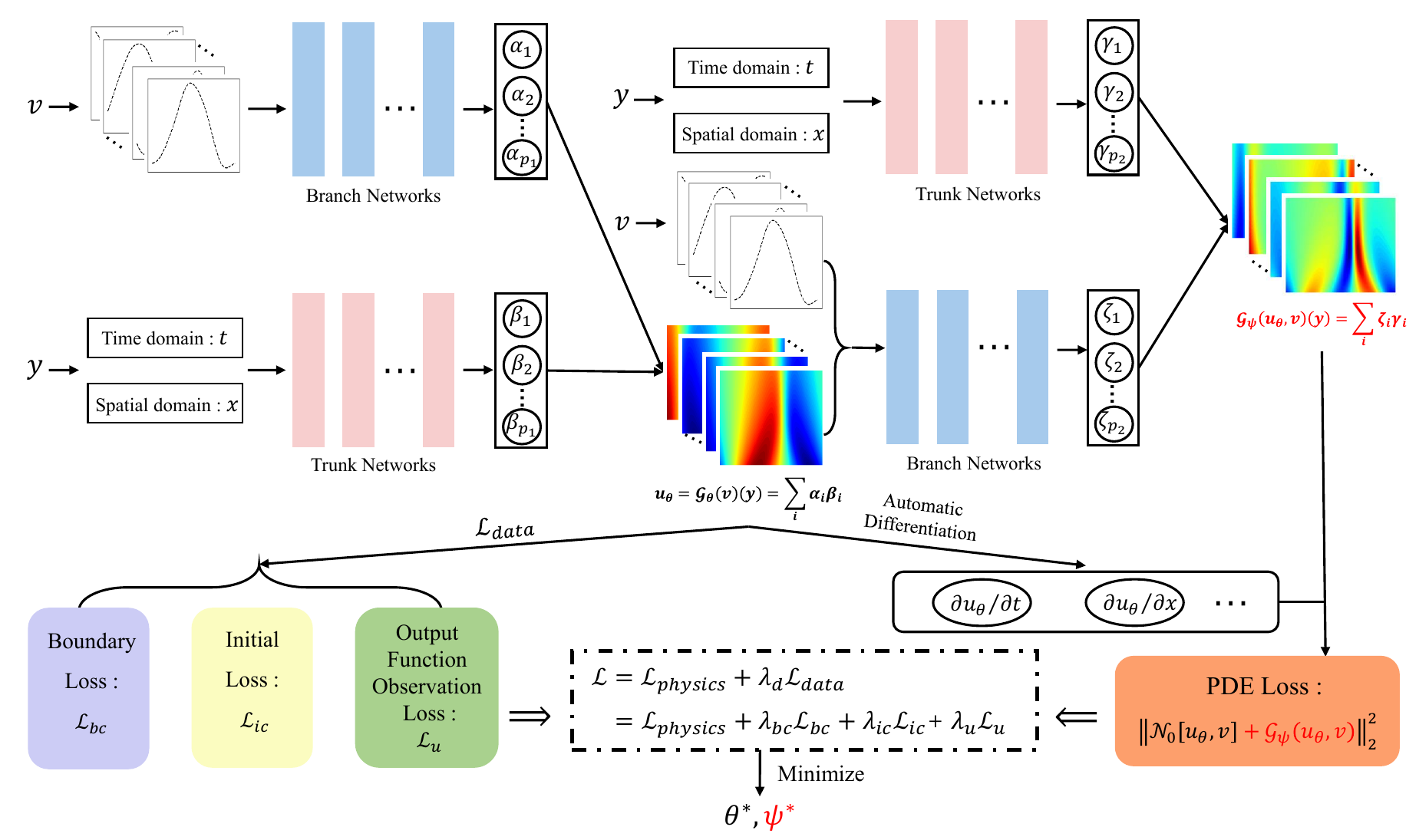}} \\
	\caption{\textbf{Schematic of the physics-guided serial DeepONet operator correction framework.} The top figure shows the serial architecture: the solution operator $\mathcal{G}_{\theta}$ maps the input function $v$ to a predicted solution $u_{\theta}(y)=\mathcal{G}_{\theta}(v)(y)$. The correction operator $\mathcal{G}_{\psi}$ then takes the observed input $v$ together with the predicted solution $u_{\theta}$ and outputs a correction term $\mathcal{G}_{\psi}(u_{\theta}, v)(y)$. The bottom figure illustrates the training objective: the data-consistency loss is imposed on available observations, while the physics loss is computed from the corrected residual $\mathcal{N}_0[u_{\theta}, v] + \mathcal{G}_{\psi}(u_{\theta}, v)$.}\label{fig:serial-deeponet}
\end{figure}

We next discuss the construction of the loss function, which consists of two parts. The first part corresponds to the physical constraint. Since $\mathcal{G}_{\theta}$ is intended to approximate the solution function, the physics loss is still built upon the prior operator $\mathcal{N}_0$, while a learned correction term is introduced to compensate for model misspecification. Specifically, the physics loss is defined by
\begin{equation} \label{eq:physics_loss}
	\mathcal{L}_{\mathrm{physics}} = \left\|
	\mathcal{N}_0[u_{\theta}(y_f), v] + \mathcal{G}_{\psi}(u_{\theta}, v)(y_f) \right\|_2^2,
\end{equation}
where $y_f$ denotes the collocation points at which the residual is evaluated. In contrast to standard physics-informed operator learning, the present formulation explicitly incorporates the output of the correction network into the physical constraint, thereby modifying the prior residual and improving the ability of the model to adapt to the true system while still retaining the available physical prior.

In addition to the physical constraint, available solution information is incorporated 
through a data consistency loss. Let $y_{\mathrm{obs}} \in \mathbb{R}^{d_y}$ denote the observation points that restricts the solution function to the available measurement locations. The data loss is defined as:
\begin{equation} \label{eq:data_loss}
	\mathcal{L}_{\mathrm{data}} =
	\left\| \mathcal{G}_{\theta}(v)(y_{\text{obs}})-u(y_{\text{obs}}) \right\|_2^2.
\end{equation}
This term may include boundary data, initial data, and interior observations, depending on 
the problem setting. The overall loss may be written in the compact form
\begin{equation*}
	\mathcal{L} = \mathcal{L}_{\text{physics}} + \lambda_{\text{d}} \mathcal{L}_{\text{data}},
\end{equation*}
where $\lambda_{\mathrm{d}}>0$ balances the physical regularization and the available 
solution data. When the available data are explicitly separated into boundary, initial, and interior 
observation subsets, the same data-consistency term can be decomposed as $
	\mathcal{L}_{\text{data}} = \lambda_{\text{bc}}\mathcal{L}_{\text{bc}} + \lambda_{\text{ic}}\mathcal{L}_{\text{ic}} + \lambda_{u}\mathcal{L}_{\text{u}}$,
which yields the equivalent refined objective
\begin{equation*}
		\mathcal{L} = \mathcal{L}_{\text{physics}} + \lambda_{\text{bc}}\mathcal{L}_{\text{bc}} + \lambda_{\text{ic}}\mathcal{L}_{\text{ic}} + \lambda_{u}\mathcal{L}_{\text{u}}.
\end{equation*}
Here, $\mathcal{L}_{\text{bc}}$, $\mathcal{L}_{\text{ic}}$, and $\mathcal{L}_{\textbf{u}}$ denote 
the losses associated with boundary conditions, initial conditions, and interior solution 
measurements, respectively, and $\lambda_{\text{bc}},\lambda_{\text{ic}},\lambda_{\textbf{u}}>0$ are the corresponding weights. Thus, the second formulation does not remove the data term; rather, it expands $\mathcal{L}_{\text{data}}$ into its constituent components. Finally, the parameters $\theta$ and $\psi$ are optimized jointly by minimizing the total loss. The corresponding schematic structure is illustrated in \cref{fig:serial-deeponet}, and the associated algorithm is summarized in \cref{alg:operator-correction}.

\begin{algorithm}[t]
	\caption{Physics-Guided Serial DeepONet for Operator Correction}
	\label{alg:operator-correction}
    \begin{algorithmic}[1]
     \STATE \textbf{Input:} Training dataset $\{(v_i, u_i)\}_{i=1}^{N_p}$, collocation points 
		$\{y_f^{(i)}\}_{i=1}^{N_f}$, observations 
		$\{(y_{\text{obs}}^{(i)},\, u(y_{\text{obs}}^{(i)}))\}_{i=1}^{N_{obs}}$, with 
		$N_{obs}=N_i+N_b+N_u$, sensor locations $\{y_c^{(i)}\}_{i=1}^{n}$, number of epochs $K$
    \FOR{$k \leftarrow 0$ to $K-1$}
    \STATE Evaluate the solution operator prediction $u_{\theta}(y) = \mathcal{G}_{\theta}(v)(y)$ and correction operator $\mathcal{G}_{\psi}(u_{\theta}(y_c), v)(y_f)$
    \STATE Compute the physics loss in \cref{eq:physics_loss}
    \STATE Compute the data-consistency loss in \cref{eq:data_loss}
    \STATE Form the total loss $\mathcal{L} = \mathcal{L}_{\text{physics}}+\lambda_{\text{d}} \mathcal{L}_{\text{data}}$ and update the networks parameters $\theta$ and $\psi$
    \ENDFOR
    \STATE \textbf{return} $\mathcal{G}_{\theta}$ and $\mathcal{G}_{\psi}$
    \end{algorithmic}
\end{algorithm}

\section{Numerical examples}\label{sec:example}
In this section, we evaluate the performance of our model on a set of benchmark problems, including a one-dimensional diffusion--reaction equation, a one-dimensional nonlinear Burgers' equation, and two two-dimensional problems: cavity flow and a hyperelasticity model. Across the main benchmark problems, we introduce representative model-form errors to assess the proposed correction framework. To further assess robustness, we investigate the effects of imperfect data by considering noisy and sparse observations in selected test cases. 

We begin by providing a detailed description of the network architectures, hyperparameter choices, and training configurations used in all numerical experiments. To ensure a fair comparison, the known model, the misspecified model, and our proposed corrected model employ networks of identical size within each example.  Moreover, both the solution operator and the correction operator are implemented using DeepONet architectures with matched capacity; in particular, their branch and trunk networks share the same depth and width. Across all experiments, we adopt the \texttt{tanh} activation function and optimize the network parameters using Adam algorithm \cite{kingma2014adam} with $\beta_1=0.999$ and $\beta_2=0.999$. The DeepONet truncation parameter is fixed to $p_1=p_2=100$ for both the branch and trunk networks. The architecture and training setup for each section are summarized in \cref{tab:nn-training-settings}. Model accuracy is quantified using the relative $L_2$ error and the mean squared error (MSE), defined as 
\begin{equation*}
		\mathrm{Relative \ L_2 \ error} = \frac{\|u - u_{\theta}\|_2}{\|u\|_2} = \frac{\left(\sum_{i=1}^{N}\left(u_i - u_{\theta, i}\right)^2\right)^{1/2}}{\left(\sum_{i=1}^{N}u_i^2\right)^{1/2}}, \quad \mathrm{MSE}
		=\frac{1}{N}\sum_{i=1}^{N}(u_i - u_{\theta, i})^2.
\end{equation*}
Unless otherwise stated, the errors reported in the tables are averaged over all samples in the test set.

\begin{table}[H]
	\centering
	%\small
	%\setlength{\tabcolsep}{10pt}
	%\renewcommand{\arraystretch}{1.8}
	\caption{Network architectures and training settings for all numerical experiments.}
	\label{tab:nn-training-settings}
	\begin{tabular}{l c c c c}
		\toprule
		Problem & Depth & Width & Learning rate & Epochs\\
		\midrule
		Diffusion-reaction & 4 & 64  & $1\times 10^{-3}$ & $1 \times 10^{5}$ \\
		Burger equation & 4 & 128 & $1\times 10^{-4}$ & $4 \times 10^{5}$ \\
		Cavity flow & 4 & 128 & $1\times 10^{-3}$ & $2 \times 10^{5}$ \\
		Hyperelasticity & 4 & 256 & $1\times 10^{-3}$ & $1 \times 10^{5}$ \\
		\bottomrule
	\end{tabular}
\end{table}

\subsection{1d diffusion-reaction equation}\label{sec:1d-dr}
We first consider the following nonlinear diffusion–reaction system presented in \cite{meng2022learning}: 
\begin{align*}
		 D \partial_x^2 u - k_r u  &= v, \quad x \in[-1,1], \\
		 u(-1) = u(1) &= 0,
\end{align*}
where $u$ denotes the solute concentration, $D=0.1$ is the diffusion coefficient, $k_r$ is the reaction rate, and $v$ is the source term. We generate a dataset using the method of manifactured solutions by sampling the reference solutions as
\begin{equation*}
	u = \frac{\left(x^2-1\right)}{10} \sum_{i=1}^5\left[\omega_{2 i-1} \sin (i \pi x)+\omega_{2 i} \cos (i \pi x)\right], \quad \omega_i\sim \mathcal{U}_{[0,1]}.
\end{equation*}
We set the reaction rate as a function of the concentration, $k_r = 0.5 \exp(-u)$, and compute the corresponding source term $v$ from the governing equation accordingly.

\begin{figure}[htbp]
	\centering	
	\subfloat[\centering : Relative $L_2$ error of $u$ and $\phi$;]{\includegraphics[width=0.4\textwidth]{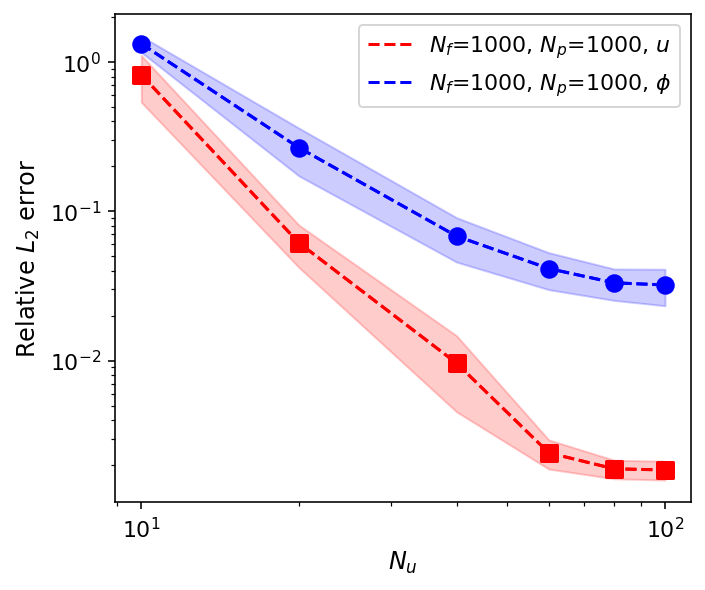}}\hspace{1.5cm}
	\subfloat[\centering : Relative $L_2$ error of $v$;]{\includegraphics[width=0.4\textwidth]{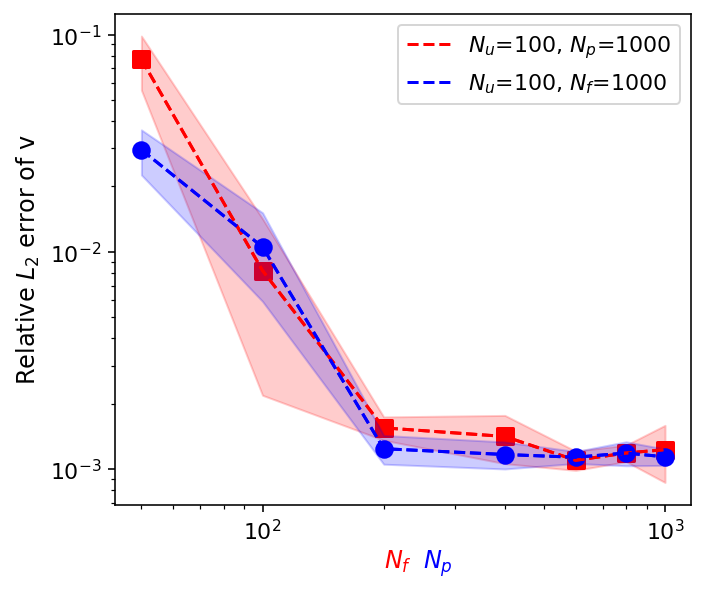}} \\
	\caption{\textbf{1d diffusion-reaction equation.} (a):Relative $L_2$ errors between corrected model predictions and the reference solution for $u$ (red) and $\phi$ (blue) as $N_u$ varies, with $N_f=1000$ and $N_p=1000$; (b): Relative $L_2$ error of the source term $v$ versus $N_p$ (blue, with $N_u=100$ and $N_f=1000$) and versus $N_f$ (red, with $N_u=100$ and $N_p=1000$). Shaded regions indicate one standard deviation over five independent experiments.}\label{fig:1d-dr-compare}
\end{figure}

In this example, we focus on learning the solution operator $\mathcal{G}: v \mapsto u$, i.e., mapping the source term to the concentration function. To demonstrate the effect of model-form error and the effectiveness of operator correction, we compare the following four settings:
\begin{itemize}
	\item \textbf{Known model:} training is carried out using the correct governing operator $\mathcal{N}$ corresponding to $D \partial_x^2 u - k_r u = v$;
	\item \textbf{Misspecified model:} training is carried out using a prior operator $\mathcal{N}_0$ given by $D \partial_x^2 u - k_r = v$, i.e., the multiplicative term $-k_r u$ is incorrectly simplified to $-k_r$;
    \item \textbf{Standard DeepONet:} solution operator trained only with the data-consistency loss;
	\item \textbf{Corrected model:} the prior operator $\mathcal{N}_0$ is augmented with a learned correction term during training.
\end{itemize}
Under our operator correction framework, the resulting corrected model is written as
\begin{equation*}
		D \frac{\partial^2 \mathcal{G}_{\theta}(v)}{\partial x^2} - k_r + \mathcal{G}_{\psi}(u_{\theta}, v) = v,
\end{equation*}
where $\mathcal{G}_{\theta}(v)$ is the solution operator surrogate for predicting $u$ from $v$, and $\mathcal{G}_{\psi}(u_{\theta}, v)$ is a correction operator intended to approximate the model discrepancy induced by misspecification. For $\mathcal{G}_{\theta}$, the branch network takes $101$ equally spaced measurements of $v$ over $[-1,1]$ as input. For the correction operator $\mathcal{G}_{\psi}$, the branch network takes the same measurements of $v$ together with the predicted solution $\mathcal{G}_{\theta}(v)$ evaluated at $101$ equally spaced locations. For error evaluation, we define the correction quantity as $\phi=-k_r u$, and all reported errors and visualizations of \(\phi\) are computed with respect to this quantity.

We investigate the effects of three factors on model performance: the number of training samples $N_p$, the number of solution observations $N_u$, and the number of collocation points $N_f$. To assess generalization, the observation locations for $u$ and the collocation points are randomly sampled from $[-1,1]$. With 1,000 training samples and 1,000 collocation points, \cref{fig:1d-dr-compare}(a) report the relative $L_2$ errors of the predicted solution $u$ and the correction term $\phi$ for operator correction model under different numbers of solution observations. The results indicate that the errors for both $u$ and $\phi$ stabilize when $N_u=100$. \cref{fig:1d-dr-compare}(b) further illustrates how the relative $L_2$ error of the source term $v$ varies with the number of collocation points and the number of training samples under different settings. We observe that $N_f=1000$ together with $N_p=1000$ is sufficient to achieve stable predictive performance. Balancing accuracy and computational cost, we adopt 1,000 training samples, 1,000 collocation points, and 100 solution observations in all subsequent examples. The test set contains 100 randomly generated samples. For each test sample, observations are provided at 201 uniformly distributed locations on $[-1,1]$.

\begin{figure}[htbp]
	\centering	
	\subfloat[\centering : Known model;]{\includegraphics[width=0.32\textwidth]{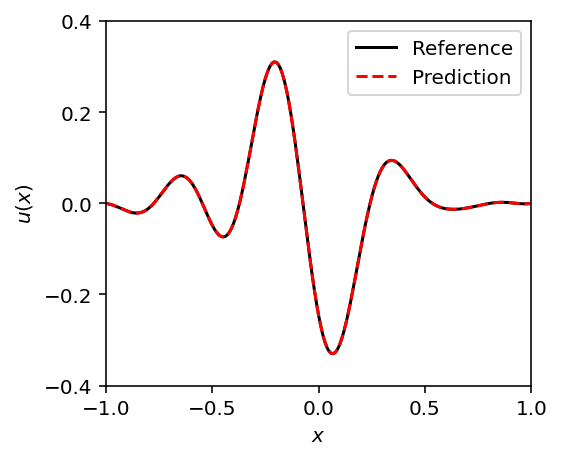}}\hspace{0.25cm}
	\subfloat[\centering : Misspecified model;]{\includegraphics[width=0.32\textwidth]{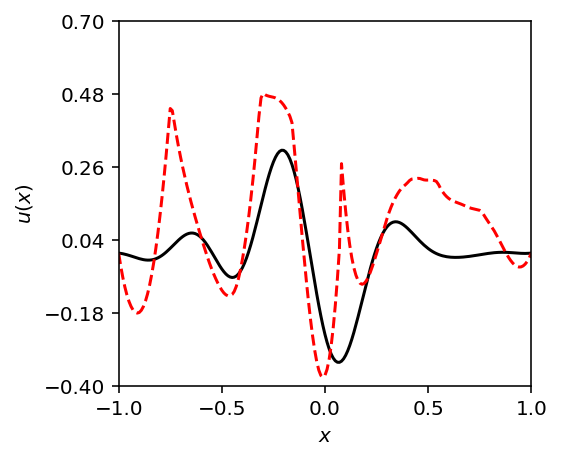}} \hspace{0.25cm}
	\subfloat[\centering : Corrected model;]{\includegraphics[width=0.32\textwidth]{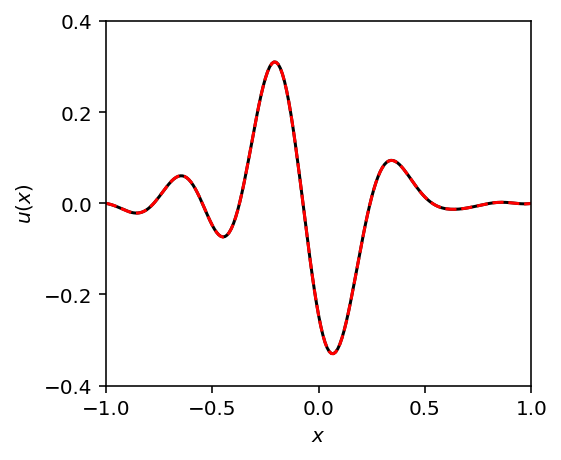}}
	\caption{\textbf{1d diffusion-reaction equation.} Reference solution $u(x)$ (black solid line) and prediction (red dashed line) for a representative test sample, (a)-(c) show the results for the known, misspecified, and corrected models, respectively.}\label{fig:1d-dr-results}
\end{figure}

In \cref{fig:1d-dr-results}, we compare the predictions of the solution $u$ for a representative test sample obtained from the known model, the misspecified model, and the corrected model, against the reference solution. Even with sufficiently dense observations of $u$ and $v$, the misspecified model exhibits significant deviations from the reference solution, indicating the inconsistency introduced by the incorrect physics constraint. In contrast, the corrected model matches the reference solution closely and substantially improves robustness by compensating for model-form errors. Quantitative results are summarized in \cref{tab:1d-dr-error}, which reports the relative $L_2$ errors of $u$, $\phi$, and $v$ for the three models. Although the corrected model does not fully reach the accuracy of the known model baseline, it significantly mitigates the impact of operator misspecification and markedly improves predictive performance. Compared with the standard DeepONet baseline, the results further show that purely data-driven operator learning can partially reduce the error induced by misspecified physics, but remains less accurate than the corrected model. This indicates that the proposed correction strategy benefits not only from learning directly from observations, but also from the explicit compensation of operator discrepancy in the governing equation.

\begin{table}[htbp]
	\centering
	\caption{\textbf{{1d diffusion-reaction equation.}} Relative $L_2$ errors between the reference and predicted solutions for $u$, $f$, and $\phi$ under the four models, where the values are reported as mean $\pm$ standard deviation over five independent runs. All reported values are scaled by $1 \times 10^{-3}$.}\label{tab:1d-dr-error}
	\begin{tabular}{ l c c c}
		\toprule
		Model & $u$ & $v$ & $\phi$\\ 
		\midrule
		Known model & $0.90 \pm 0.14$ & $0.58 \pm 0.04$ & $0.81 \pm 0.13$ \\ 
		Misspecified model & $1614.20 \pm 245.999$ & $109386.97 \pm 155420.55$ & $1209.27 \pm 52.66$ \\ 
        Standard DeepONet & $3.60 \pm 0.89$ & --- & --- \\ 
		Corrected model & $1.85 \pm 0.27$ & $1.10 \pm 0.13$ & $32.10 \pm 8.81$ \\ 
		\bottomrule
	\end{tabular}
\end{table}

\begin{figure}[t]
	\centering	
	\subfloat[\centering : Sparse observations;]{\includegraphics[width=0.4\textwidth]{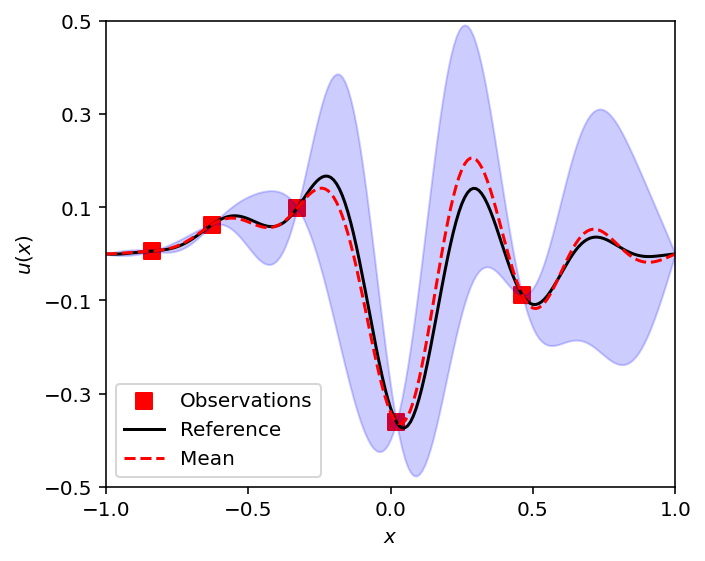}}\hspace{1.0cm}
	\subfloat[\centering : Noisy observations;]{\includegraphics[width=0.4\textwidth]{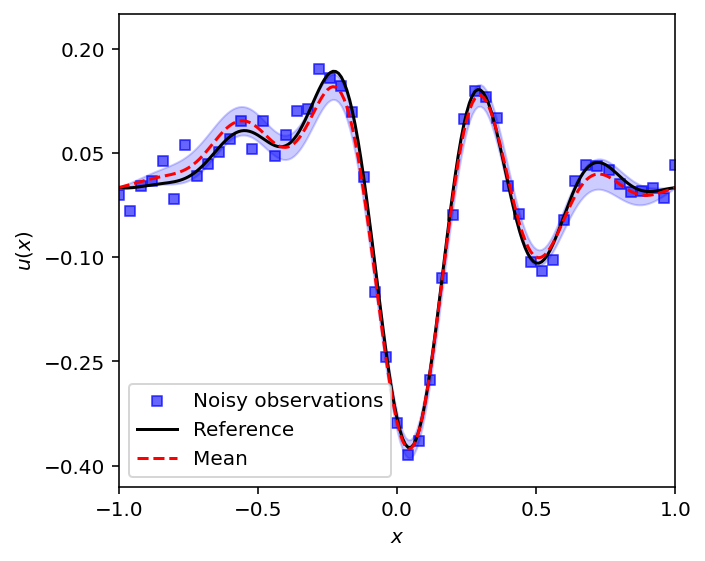}} \\
	\caption{\textbf{1d diffusion-reaction equation under imperfect observations.} (a): Prediction of $u(x)$ with five sparse observations; (b): Prediction of $u(x)$ with noisy observations at noise level $\sigma=0.02$. The black solid line denotes the reference solution, the red dashed lines denote the ensemble mean prediction, and the shaded regions indicate two standard deviations estimated from ten independent runs. The square markers indicate the available observations (red squares for sparse observations and blue squares for noisy observations).}\label{fig:1d-dr-grappy-noisy}
\end{figure}

In the previous experiments, we assumed that observations were sufficiently informative for accurate evaluation. We now consider more challenging scenarios to assess robustness under imperfect data: (i) sparse observations and (ii) noisy observations. For the sparse observations case, we assume that $u$ is observed only at five randomly selected spatial locations, while all other settings remain the same as in the previous experiments. For the noisy observations case, we corrupt both $u$ and $v$ with i.i.d. Gaussian noise: $u(x) = u(x) + \epsilon_u$ and $v(x) = v(x) + \epsilon_v$, where $\epsilon_u,\epsilon_v\sim\mathcal{N}(0,\sigma^2)$. We set the noise level to $\sigma=0.02$ and examine its effect on model accuracy. The results are summarized in \cref{fig:1d-dr-grappy-noisy}, where (a) corresponds to sparse observations and (b) corresponds to noisy observations. Under sparse observations, the limited measurement information leads to increased predictive uncertainty in the solution $u$. Nevertheless, the proposed operator correction model yields a mean prediction that remains in good agreement with the reference solution, suggesting that the learned correction can still provide useful predictive information under limited data. Under noisy measurements, the model remains stable and accurately captures the main profile of u, demonstrating its robustness to observation noise. The uncertainty estimates are obtained using deep ensembles, constructed from ten independently trained models with different random seeds.

\subsection{1d burgers' equation}\label{sec:1d-bu}
Burgers’ equation is a canonical nonlinear convection–diffusion PDE and has been widely used as a simplified model in areas such as gas dynamics and turbulence \cite{burgers2013nonlinear, whitham2011linear}. The nonlinear advection term $u\,\partial_x u$ and the viscous dissipation term $\nu\,\partial_x^2 u$ jointly shape the solution dynamics, including the formation and smoothing of sharp gradients \cite{leveque2002finite, bec2007burgers}. Owing to its well-known nonlinear behaviors and the availability of reliable numerical solvers, it has become a standard benchmark for operator learning and physics-informed modeling. In this subsection, following the benchmark setting in \cite{wang2021learning}, we consider the one-dimensional nonlinear Burgers’ equation:
\begin{equation}\label{eq:2-1d-burgers}
	\begin{aligned}
		\frac{\partial u}{\partial t}+u \frac{\partial u}{\partial x} & = \nu \frac{\partial^2 u}{\partial x^2}, & & (x,t)\in(0,1)\times(0,1],\\
		u(x,0) & = v(x),  && x\in(0,1),\\
		\end{aligned}
\end{equation}
with periodic boundary conditions: $u(0,t)  = u(1,t)$ and $\frac{\partial u}{\partial x}(0,t)  = \frac{\partial u}{\partial x}(1,t)$. Here, $\nu$ denotes the viscosity coefficient, which is set to $\nu=0.01$. Our goal is to learn the solution operator that maps the initial condition $v$ to the corresponding spatio-temporal solution function $u(x,t)$ over the entire domain, i.e., $\mathcal{G}: v \mapsto u$. We note that, in this example, the input function $v$ for the solution operator is a temporal slice of the solution function $u(x,t)$ at $t=0$ (i.e., the initial condition). Nevertheless, to maintain a unified formulation of the proposed framework, we still define the correction operator to take $(v, u)$ as its inputs. 

\begin{figure}[t]
	\centering	
	{\includegraphics[width=0.4\textwidth]{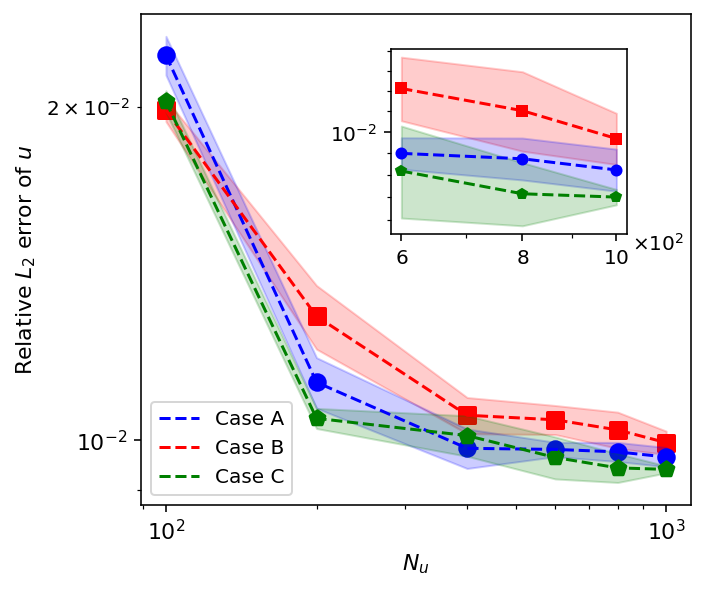}}
	\caption{\textbf{1d burgers' equation.} Relative $L_2$ error of the solution $u$ versus the number of solution observations $N_u$ for Cases A–C under the corrected model. The results are obtained with $N_p=2000$ training samples and an $N_f=101 \times 101$ collocation grid. Markers and shaded regions denote the mean and one standard deviation, respectively, computed over five independent runs.}\label{fig:1d-bugers-compare}
\end{figure}

To evaluate the proposed framework under different types of model-form errors, we consider three misspecification scenarios:
\begin{itemize}
	\item \textbf{Case A (perturbed advection nonlinearity):} The governing equation is assumed to contain an additional cubic term $\epsilon u^3$ \cite{ebers2024discrepancy}:
	\begin{equation*}
		\frac{\partial u}{\partial t}+u \frac{\partial u}{\partial x} + \epsilon u^3 = \nu \frac{\partial^2 u}{\partial x^2},
	\end{equation*}
	where $\epsilon=10$.

	\item \textbf{Case B (advection term omitted):} The nonlinear advection term is neglected, resulting in a purely diffusive model:
	\begin{equation*}
		\frac{\partial u}{\partial t} = \nu \frac{\partial^2 u}{\partial x^2},
	\end{equation*}
	consistent with the setup in \cite{podina2022pinn}.
	
	\item \textbf{Case C (diffusion term omitted):} The viscous dissipation term \(\nu\,\partial_x^2 u\) is removed, yielding the inviscid Burgers’ equation:
	\begin{equation*}
		\frac{\partial u}{\partial t}+u \frac{\partial u}{\partial x} = 0.
	\end{equation*}
\end{itemize}
In our proposed model, the correction operator $\mathcal{G}_{\psi}$ is appended to the part of the governing equation affected by model misspecification. Under this convention, the corresponding target correction terms for cases A–C are, respectively,
$-\epsilon u^3$, $u\,\frac{\partial u}{\partial x}$, and $\nu\,\frac{\partial^2 u}{\partial x^2}$.

The initial condition $v$ is drawn from a Gaussian random field as $v\sim\mathcal{N}(0,\,25^2(-\Delta+5^2 I)^{-4})$. Reference solutions $u(t,x)$ are obtained by integrating the governing equation using an explicit forward Euler scheme in time together with a Fourier pseudo-spectral discretization in space. We construct a dataset comprising $2100$ input–output function pairs, of which $2000$ are used for training and the remaining $100$ are reserved for testing. Each realization is evaluated on a uniform grid of $101$ time instances and $201$ spatial locations. For physics constraints, the PDE residual is evaluated on a uniform $101\times 101$ spatio-temporal collocation grid over $(0,1)\times(0,1]$. We take $101$ equally spaced observations of $v$ on $[0,1]$, which serve both as the branch input of the solution operator and the initial training set, i.e., $N_i=101$. For the correction operator, the solution prediction $u_{\theta}$ is incorporated into the branch input through its evaluations on a fixed $51\times 51$ spatio-temporal grid. Periodic boundary constraints are imposed at $x=0$ and $x=1$ using $100$ temporal training points sampled over the interval $[0,1]$. The loss weights are fixed as $\lambda_{\text{bc}}=1$ and $\lambda_{\text{ic}}=\lambda_{\text{u}}=50$. We first investigate the dependence of the corrected model accuracy on the number of solution observations $N_u$. As shown in \cref{fig:1d-bugers-compare}, for all three cases the relative $L_2$ error of the solution stabilizes at $N_u=800$. Accordingly, we set $N_u=1000$ in the subsequent experiments.

\begin{figure}[t]
	\centering	
	{\includegraphics[height=0.53\textheight, width=1.0\textwidth]{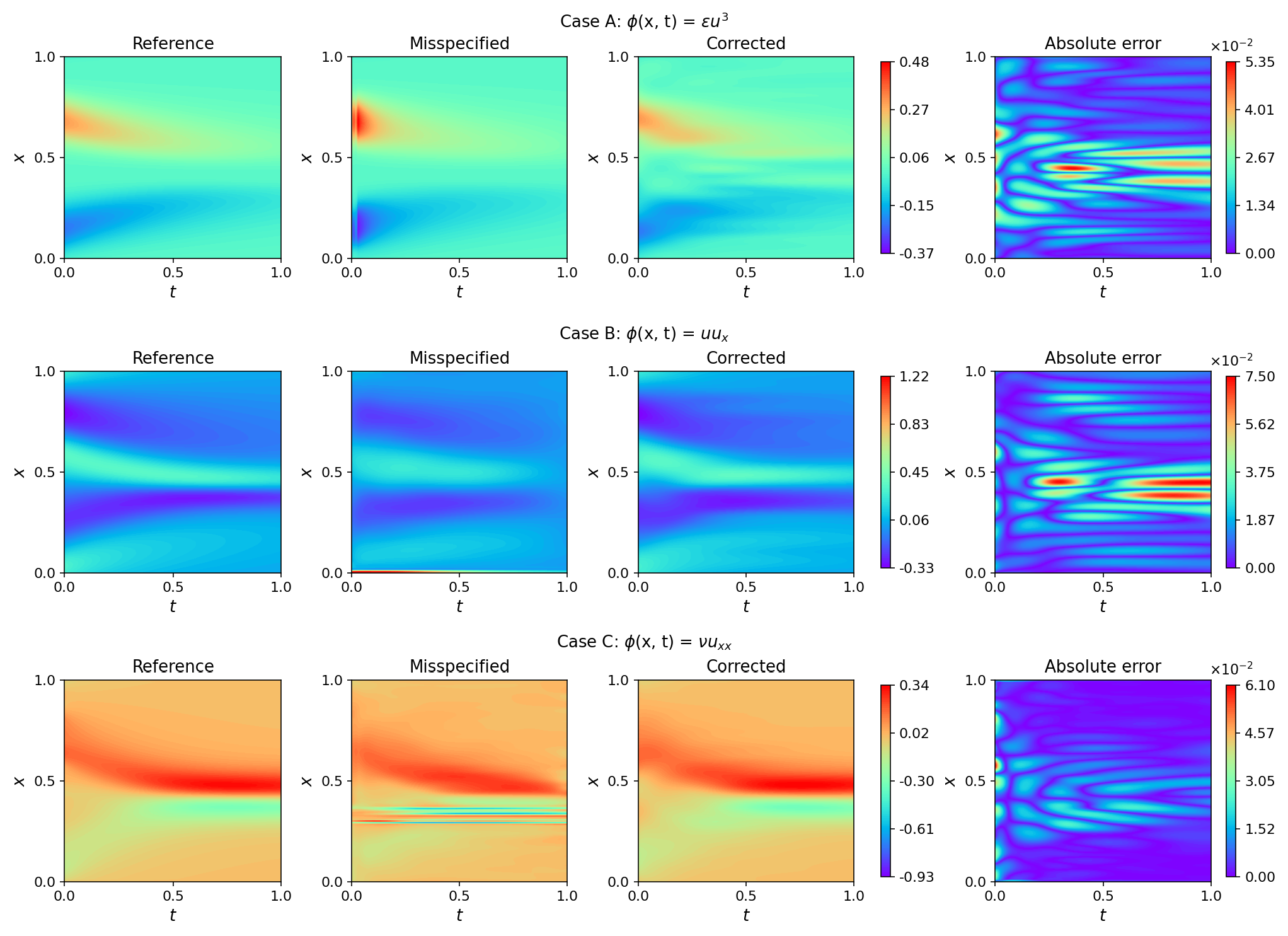}}
	\caption{\textbf{1d burgers' equation.} Correction term $\phi(x,t)$ for a representative test sample in three cases. Columns 1–3 display the reference solution and the corresponding predictions of the misspecified and corrected models, respectively; column 4 reports the absolute error of the corrected model relative to the reference. }\label{fig:1d-bugers-phi}
\end{figure}

\cref{fig:1d-bugers-phi} presents the numerical results of the correction term $\phi$ for a representative test sample under the three misspecification cases. We compare the reference correction term with the predictions obtained from the misspecified and corrected models, where the model settings are defined in \cref{sec:1d-dr}. When the governing equation is incorrectly specified, the misspecified model produces correction terms that deviate substantially from the reference in both magnitude and spatio-temporal distribution. In contrast, the corrected model provides much more accurate approximations of $\phi$ across all three cases, demonstrating its ability to learn the discrepancy induced by model misspecification.

\begin{table}[htbp]
	\centering
	\caption{\textbf{1d burgers' equation.} Relative $L_2$ errors between the reference and model predictions for $u$ and $\phi$ under the known model, standard DeepONet, misspecified model, and corrected model in Cases A--C. Values are reported as mean $\pm$ standard deviation over five independent runs. The errors for $u$ are expressed in percentages, whereas the errors for $\phi$ are reported as raw relative $L_2$ errors.}
	\label{tab:1d-burgers-error}
	\begin{tabular}{llcc}
		\toprule
		\multicolumn{2}{c}{Model} & $u$ error ($\%$) & $\phi$ error \\
		\midrule
		\multicolumn{2}{c}{Known model}
		& $0.79 \pm 0.06$
		& --- \\
		\midrule
        \multicolumn{2}{c}{Standard DeepONet}
		& $2.07 \pm 0.50$
		& --- \\
		\midrule
		
		\multirow{2}{*}{Case A}
		& Misspecified
		& $14.23 \pm 0.24$
		& $0.67 \pm 0.06$ \\
		& Corrected
		& $0.99 \pm 0.02$
		& $0.25 \pm 0.00$ \\
		\addlinespace
		
		\multirow{2}{*}{Case B}
		& Misspecified
		& $7.98 \pm 0.25$
		& $1.47 \pm 0.99$ \\
		& Corrected
		& $1.01 \pm 0.05$
		& $0.23 \pm 0.01$ \\
		\addlinespace
		
		\multirow{2}{*}{Case C}
		& Misspecified
		& $6.64 \pm 0.57$
		& $1.13 \pm 0.24$ \\
		& Corrected
		& $0.95 \pm 0.03$
		& $0.19 \pm 0.01$ \\
		\bottomrule
	\end{tabular}
\end{table}

The relative $L_2$ errors for $u$ and $\phi$ are summarized in \cref{tab:1d-burgers-error} for all three model settings across Cases A–C. When the misspecification does not explicitly alter the differential operators acting on the solution, the uncorrected baseline remains capable of reproducing the solution at a qualitative level. By contrast, omission of essential derivative terms (e.g. the advection term $u\,\partial_x u$; the diffusion term $\nu\,\partial_x^2 u$) generally precludes accurate recovery of the spatio-temporal evolution. Taken together, these results demonstrate the effectiveness of the proposed model in mitigating model-form errors by fully exploiting available observations, particularly those attributable to incorrect differential operator.

\begin{figure}[htbp]
	\centering	
	{\includegraphics[height=0.20\textheight, width=1.0\textwidth]{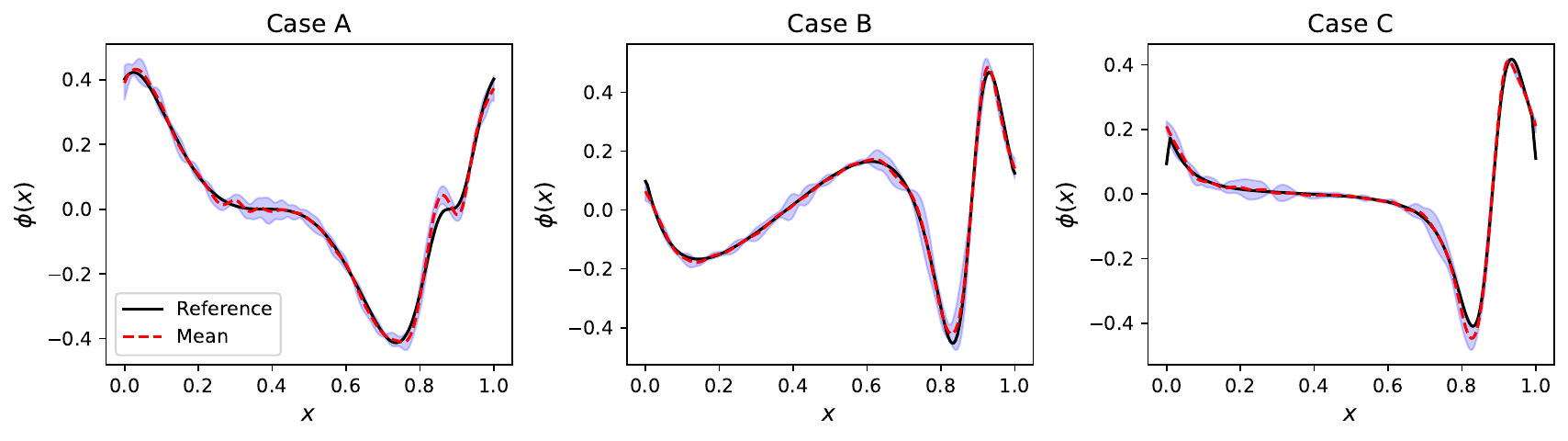}}\\
	{\includegraphics[height=0.20\textheight, width=1.0\textwidth]{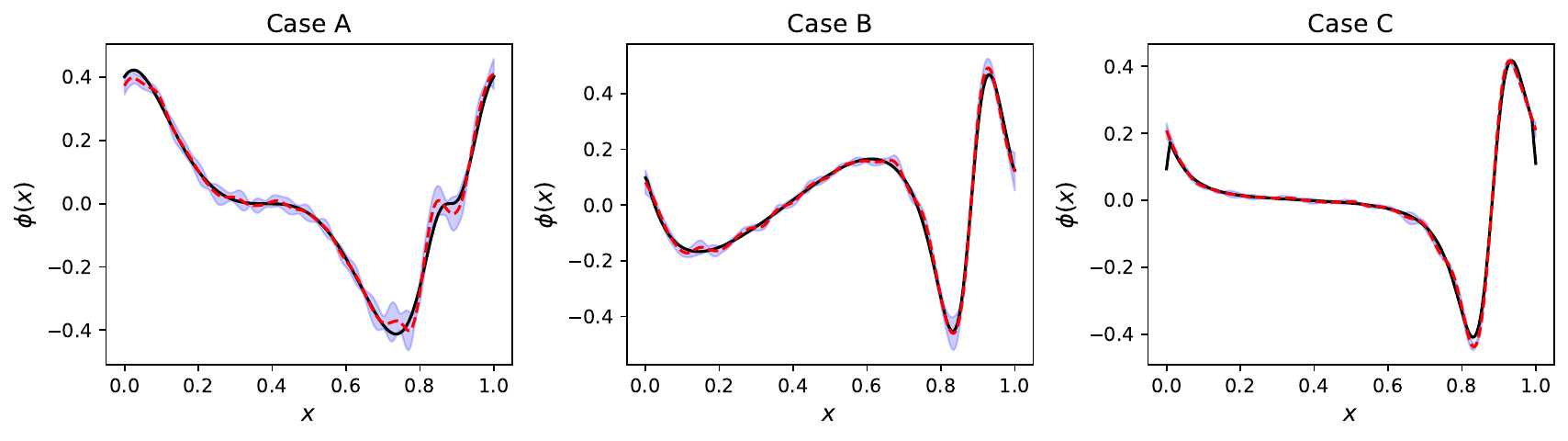}}
	\caption{\textbf{1d burgers' equation.} Correction term $\phi$ at $t=0.5$ under sparse solution observations for Cases A–C. Results are shown for $100$ (top) and $200$ (bottom) solution observations. The black line denotes the reference solution, the red dashed line denotes the ensemble mean prediction, and the shaded regions indicates two standard deviation computed from five independent runs.}\label{fig:1d-burgers-grappy}
\end{figure}

We further investigate the model performance under sparse observations and its uncertainty quantification capability. \cref{fig:1d-burgers-grappy} shows the correction term $\phi$ at $t=0.5$ for a representative test realization, obtained by training with $100$ or $200$ solution observations. For $N_u=100$, the ensemble mean displays modest localized oscillations, nevertheless, the mean prediction remains in good agreement with the reference, and the ensemble uncertainty band covers the main variations of the reference solution. When the number of observations increases to $N_u=200$, the predictive variance decreases markedly, the ensemble mean becomes smoother, and the agreement with the reference improves, indicating that additional observations effectively reduce uncertainty and enhance predictive stability. These results indicate that the proposed corrected framework remains robustness under sparse measurements and provides meaningful uncertainty quantification.

\subsection{2d cavity flow}\label{sec:2d-cav}
Newtonian and non-Newtonian fluids are fundamentally distinguished by their constitutive relation between the deviatoric stress and the rate of strain tensor. For Newtonian fluids, the viscous stress depends linearly on the velocity gradient with a constant dynamic viscosity, whereas for non-Newtonian fluids the apparent viscosity typically varies with the shear rate and may exhibit shear-thinning/thickening and other nonlinear rheological effects. Consequently, modeling a non-Newtonian fluid as Newtonian introduces systematic bias in the constitutive (viscous) term of the governing equations. Motivated by this observation, we consider a two-dimensional lid-driven cavity flow of a power-law non-Newtonian fluid.

\paragraph{Ground truth}
Let $\boldsymbol{u}=(u_x,u_y)$ denote the dimensionless velocity field and $P$ represents the dimensionless pressure on the unit square $\Omega=[0,1]\times[0,1]$. The flow satisfies the incompressible steady 
momentum balance
\begin{equation}
	\begin{aligned}
    \boldsymbol{u}\cdot\nabla\boldsymbol{u}
		&= -\nabla P
		+ \frac{1}{Re}\nabla\cdot\left(\nu \boldsymbol{S}\right),\\
		\nabla \cdot \boldsymbol{u} &= 0,
	\end{aligned}
	\label{eq:2-2d-cavity}
\end{equation}
where $\boldsymbol{S}= \nabla\boldsymbol{u}+(\nabla\boldsymbol{u})^\top $ is the symmetric velocity-gradient tensor, and $Re$ is the generalized
Reynolds number for the power-law fluid. The associated shear-rate magnitude is defined as $|\boldsymbol{S}|=\sqrt{\boldsymbol{S}:\boldsymbol{S}/2}$,
and the dimensionless power-law viscosity is modeled by $\nu = |\boldsymbol{S}|^{n-1}$, and we select $n=1.5$. For this example, the Reynolds number $Re$ is sampled from $[100,200]$. The top boundary is driven by the prescribed velocity profile
\begin{equation}
	\boldsymbol{u}(x,1) =
	\begin{bmatrix}
		v(x)\\
		0
	\end{bmatrix},
	\qquad
	v(x) = 1-
	\frac{
		\cosh\left[10\left(x-\frac{1}{2}\right)\right]
	}{
		\cosh(5)
	},
	\qquad 0\leq x\leq 1.
	\label{eq:cavity-top-bc}
\end{equation}
No-slip boundary conditions are imposed on the remaining walls, namely the bottom boundary and the two vertical side boundaries. We focus on learning the parametric solution operator that maps the Reynolds number to the velocity field $
\mathcal{G}:Re\mapsto \boldsymbol{u}$.
The training and test data are generated using a lattice Boltzmann equation (LBE) solver following Wang et al. \cite{wang2015localized}. To determine whether the iterative solver has reached a steady state, we use a relative-change criterion based on the velocity field. Let $\boldsymbol{u}_T$ 
denote the velocity field at iteration $T$, and define $
	E=
		\sum |\boldsymbol{u}_{T+500}-\boldsymbol{u}_T|
	/
		\sum |\boldsymbol{u}_T|
	,
$ where the summations are taken over all grid points and velocity components. The iteration is terminated once $E\leq 10^{-6}$, and the corresponding solution is regarded as the steady-state velocity field. To avoid excessive computational cost in rare nonconvergent cases, a maximum of $5\times10^5$ iterations is also imposed. The resulting dataset consists of $1000$ training samples and $100$ test samples, with all realizations simulated on a uniform $101\times101$ grid. 

\paragraph{Misspecified model}
At the modeling stage, we assume that the system is incorrectly treated as Newtonian fluid. That is, the shear rate dependent viscosity is replaced by a constant viscosity. Using the same nondimensional formulation as \cref{eq:2-2d-cavity}, the resulting misspecified Newtonian-form momentum equation is
\begin{equation}
	\boldsymbol{u}\cdot\nabla\boldsymbol{u} = -\nabla P + \frac{1}{Re}\nabla^2\boldsymbol{u}.
	\label{eq:newtonian-misspecified-model}
\end{equation}
In the original dimensional formulation used before nondimensionalization, the constant viscosity coefficient in the misspecified Newtonian model is taken as $\nu_0=10^{-4}$. After nondimensionalization, this dimensional coefficient is 
absorbed into the dimensionless scaling and the Reynolds number. Therefore, it does not appear as an additional independent coefficient in \cref{eq:newtonian-misspecified-model}. To isolate the impact of constitutive misspecification on the velocity correction, we additionally learn a pressure operator $\mathcal{G}_{\xi}:Re\mapsto P$ and provide sufficiently dense supervision for $P$, so that pressure errors do not dominate the velocity correction. The correction operator $\mathcal{G}_{\psi}(Re,\boldsymbol{u})$ is vector-valued and accounts for both velocity components. The predicted velocity field supplied to the correction operator is represented by $\mathcal{G}_\theta(v)$ evaluated on a fixed $51\times51$ uniform grid in $\Omega$. Physics constraints are enforced on a $101\times101$ uniform collocation grid. For the boundary training set, in addition to the $101$ observations on the top, we use $51$ equispaced observations on each of the remaining boundaries. Interior velocity observations consist of $250$ points, while pressure observations are provided on a $81\times81$ grid to stabilize pressure learning and disentangle its effect from the velocity correction. 

\begin{figure}[!htbp]
	\centering	
	{\includegraphics[height=0.38\textheight, width=1.0\textwidth]{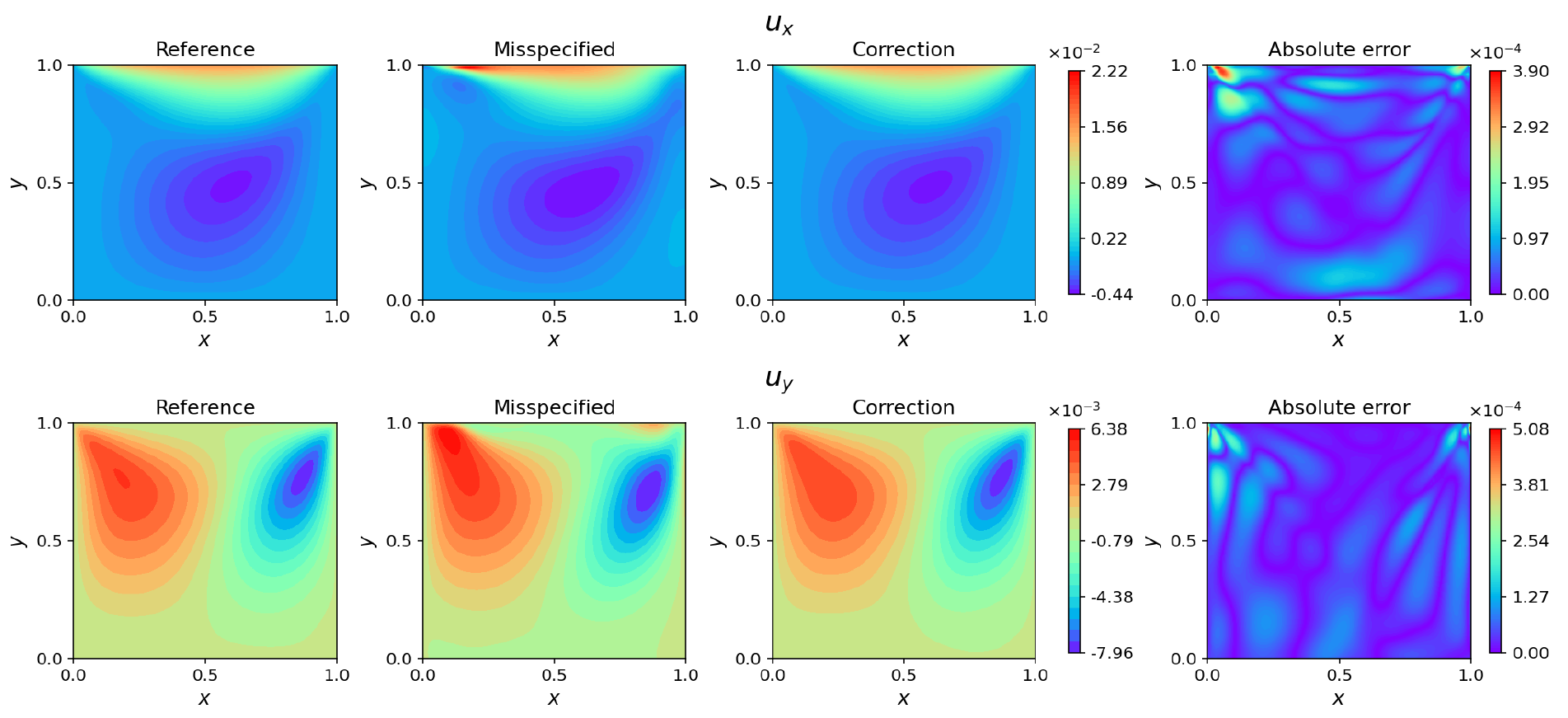}}
	\caption{\textbf{2d cavity flow.} Representative velocity field results for a test realization, where the top and bottom rows show the $x$- and $y$- components $u_x$ and $u_y$, respectively. From left to right, the columns report the reference solution, the misspecified prediction, the corrected prediction, and the absolute error of the corrected prediction relative to the reference.}\label{fig:2d-cavity-uv}
\end{figure}

\cref{fig:2d-cavity-uv} shows the contour plots of the velocity components for a representative test realization. The misspecified Newtonian model exhibits noticeable discrepancies from the reference solution, reflecting the systematic bias caused by the inaccurate constitutive assumption. By contrast, the corrected model provides velocity predictions that are in closer agreement with the reference solution and shows reduced local errors. The quantitative results are summarized in \cref{tab:2d-cavity-error}, including the relative $L_2$ errors of the velocity components $u_x$ and $u_y$ and the pressure $P$, as well as the mean squared errors of the residual components $f_x$ and $f_y$. These numerical errors are consistent with the visual observations and further confirm that the proposed operator correction framework improves the predictive accuracy under model misspecification.

\begin{table}[H]
\centering
\caption{\textbf{2d cavity flow.} Relative $L_2$ errors for the velocity components $u_x$, $u_y$, and pressure $P$, together with the mean squared errors of the residual components $f_x$ and $f_y$, for the four models considered above. Values are reported as mean $\pm$ standard deviation over five independent runs and are expressed in percentages.}\label{tab:2d-cavity-error}

\begin{tabular}{lccccc}
	\toprule
	Model & $u_x$ & $u_y$ & $P$ & $f_x$ & $f_y$ \\
	\midrule
    
    Known model
	& $0.66 \pm 0.04$
	& $0.93 \pm 0.07$
	& $0.13\% \pm 0.01\%$
	& $0.18\%  \pm 0.43\%$
	& $0.19\%  \pm 0.25\%$ \\
	
	Misspecified
	& $22.09 \pm 0.48$
	& $26.58 \pm 0.83$
	& $0.01 \pm 0.01\%$
	& $0.41 \pm 0.04$
	& $0.55 \pm 0.04$ \\
	%\addlinespace

    Standard DeepONet
	& $1.77 \pm 0.48$
	& $1.50 \pm 0.21$
	& ---
	& ---
	& --- \\
	
	Corrected
	& $0.82 \pm 0.23$
	& $1.11 \pm 0.53$
	& $0.80\% \pm 0.04\%$
	& $0.77\% \pm 0.38\%$
	& $0.48\% \pm 0.21\%$ \\
	\bottomrule
\end{tabular}
\end{table}

\subsection{2d hyperelastic problem}\label{sec:2d-hyperelastic}
Hyperelastic material models are commonly employed to describe materials that undergo large deformations while still exhibiting a recoverable elastic response. Their constitutive behavior is typically characterized by a strain-energy density function, which makes them well suited for capturing the nonlinear mechanical behavior of rubber-like materials, elastomers, and certain soft biological tissues \cite{ogden1997non, holzapfel2002nonlinear}. Finally, to further assess the effectiveness of the proposed operator correction model in a nonlinear elasticity setting, we consider a two-dimensional hyperelastic beam under compression \cite{bouziani2024structure}. Specifically, the displacement field $\boldsymbol{u}=(u_x,u_y)$ is governed by the nonlinear elasticity equation:
\begin{equation*}
		-\nabla \cdot \boldsymbol{P}(\boldsymbol{u}) = \boldsymbol{f}, \qquad \text{in } \Omega,
\end{equation*}
where the body force is taken as $\boldsymbol{f}=(0,-1000)^{\top}$.
The first Piola-Kirchhoff stress tensor is defined by
\begin{equation*}
		\boldsymbol{P} = \mu \boldsymbol{F} + (-\mu + \lambda\ln J) \boldsymbol{F} ^{- \mathsf T}, \quad \boldsymbol{F} = \boldsymbol{I}+\nabla \boldsymbol{u}, 
		\qquad
		J=\det(\boldsymbol{F}).
\end{equation*}
The Lam\'e parameters are determined by the Young's modulus $E=10^6$ and Poisson ratio $\nu=0.3$ through
\[
\mu=\frac{E}{2(1+\nu)},
\qquad
\lambda=\frac{E\nu}{(1+\nu)(1-2\nu)}.
\]

The computational domain is a slender rectangle $\Omega = [0,1]\times[0,0.1]$. Boundary conditions are prescribed as follows. On the left boundary $\{0\}\times[0,0.1]$, the beam is clamped and $
\boldsymbol{u}=\boldsymbol{0}$. On the top and bottom boundaries, natural traction-free conditions are imposed as $
\boldsymbol{P}(\boldsymbol{u})\boldsymbol{n}=\boldsymbol{0}$.
On the right boundary $\Gamma_R$, we apply a horizontal compressive displacement
\[
\boldsymbol{u}|_{\Gamma_R}=(v, 0)^{\mathsf T}, \quad \Gamma_R=\{1\}\times[0,0.1],
\]
where $v$ is a scalar boundary input. For simplicity, we consider constant compressive loads of the form $v=-\epsilon, \epsilon\in \mathbb{R}_{+}$. In this example, we study the mapping from a prescribed compressive load on the right boundary to the displacement field over the entire domain:
$\mathcal{G}: v \mapsto \boldsymbol{u}$.

To mitigate the adverse effects of the scale disparity among different physical quantities on network training and to improve the numerical stability in solving the hyperelastic problem, we first nondimensionalize the governing equations \cite{ren2024mixed,langtangen2016scaling}. It should be emphasized that this nondimensionalization does not involve any rescaling of the spatial coordinates; rather, it is applied only to the displacement, stress, energy density, body force, and the associated material parameters. Building on this, in order to further improve the numerical stability of network training while simultaneously exploiting the variational structure and the local equilibrium constraints of the hyperelastic problem, we adopt a joint loss construction that combines the strong-form residual with the total potential energy \cite{abueidda2023enhanced}. Specifically, the total potential energy functional of the hyperelastic body under the force term $\boldsymbol{f}$ is given by
\begin{equation*}
	\Pi(\boldsymbol{u}) =
	\int_{\Omega} W(\boldsymbol{F})\,\mathrm{d}\Omega
	-
	\int_{\Omega}\boldsymbol{f}\cdot\boldsymbol{u}\,\mathrm{d}\Omega
	-
	\int_{\Gamma_N}\bar{\boldsymbol{t}}\cdot\boldsymbol{u}\,\mathrm{d}\Gamma,
\end{equation*}
where $\bar{\boldsymbol{t}}$ denotes the prescribed surface traction on the Neumann boundary $\Gamma_N$, and $W(\boldsymbol{F}_{\theta})$ is the strain energy density function, which is given by
$W(\boldsymbol{F}) = \frac{\mu}{2}\Big(\operatorname{tr}(\boldsymbol{F}^{\mathsf T}\boldsymbol{F})-2-2\ln J\Big) +\frac{\lambda}{2}(\ln J)^2$. Since traction-free boundary conditions are imposed on the top and bottom boundaries, the total potential energy functional reduces to $\Pi(\boldsymbol{u}) = \int_{\Omega} W(\boldsymbol{F})\,\mathrm{d}\Omega - \int_{\Omega}\boldsymbol{f}\cdot\boldsymbol{u}\,\mathrm{d}\Omega$, and the energy-based loss term  can be defined by $\mathcal{L}_{\text{energy}} = \Pi(\boldsymbol{u}_{\theta})$. Accordingly, the overall loss function of the governing equation can be written as:
\begin{equation}
	\begin{aligned}
		\mathcal{L} = \mathcal{L}_{\text{physics}} + \lambda_{{\text{d}}} \mathcal{L}_{\text{data}} + \lambda_{{\text{e}}} \mathcal{L}_{\text{energy}},
	\end{aligned}
\end{equation}
where $\lambda_e$ denotes the weight coefficient associated with the energy loss term.

\begin{figure}[htbp]
	\centering	
	\subfloat[$u_x$;]{\includegraphics[height=0.36\textheight, width=1.0\textwidth]{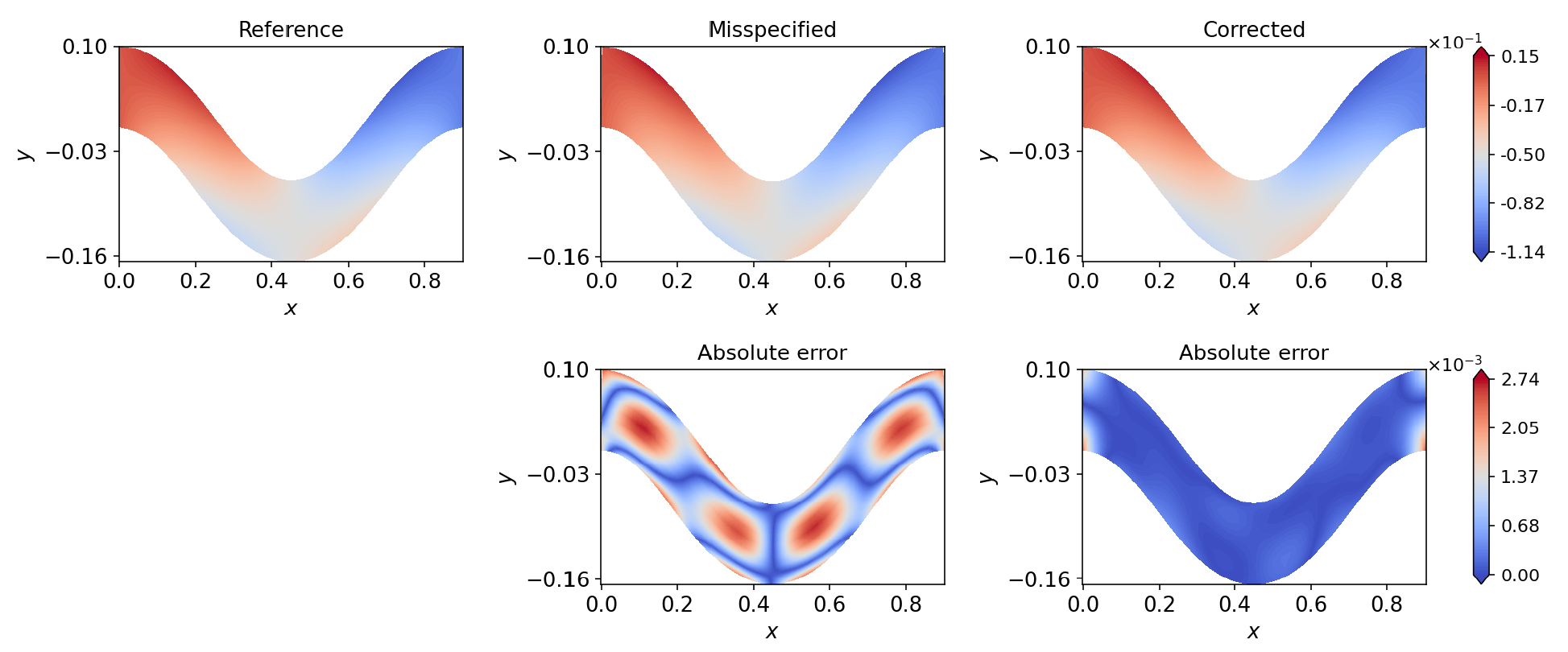}} \\
	\subfloat[$u_y$;]{\includegraphics[height=0.36\textheight, width=1.0\textwidth]{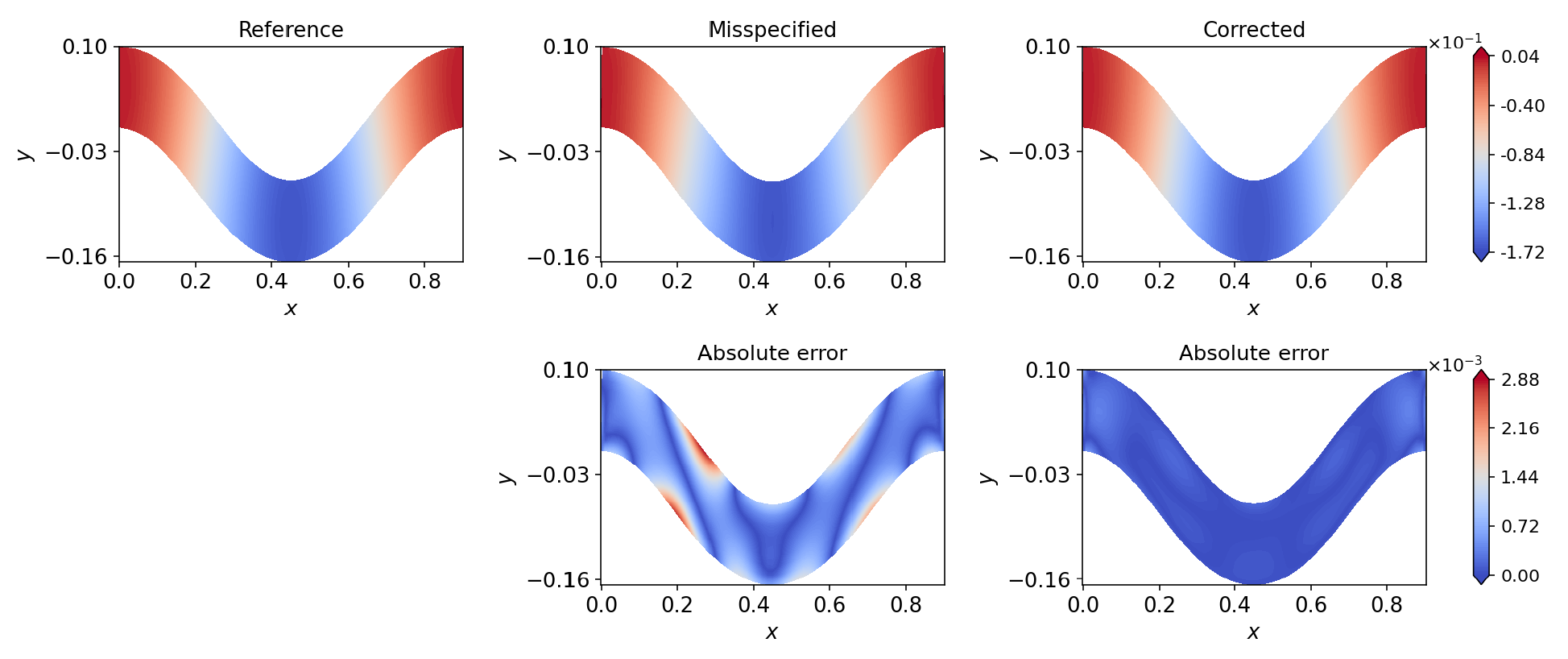}}
	\caption{\textbf{2d hyperelastic problem.} Displacement field components for a representative test sample with $\epsilon = 0.0991$: (a) $u_x$ and (b) $u_y$. For each component, the first row presents the reference solution and the predictions of the misspecified and corrected models, while the second row shows the corresponding absolute errors relative to the reference solution.}\label{fig:2d-hyper-u}
\end{figure}

Here, we consider the case where the governing equation is incorrectly specified as a linear elastic model:
\begin{equation}\label{eq:2d-linear}
	\begin{aligned}
		\begin{cases}
			\nabla \cdot \boldsymbol{\sigma}(\boldsymbol{u}) + \boldsymbol{f} = \boldsymbol{0}, & \text{in} \Omega, \\
			\boldsymbol{\sigma}(\boldsymbol{u}) = \lambda \, \mathrm{tr} (\boldsymbol{\varepsilon}(\boldsymbol{u})) \,\boldsymbol{I} + 2\mu\,\boldsymbol{\varepsilon}(\boldsymbol{u}), & \text{in} \Omega,\\
			\boldsymbol{\varepsilon}(\boldsymbol{u})=\dfrac{1}{2}(\nabla \mathbf{u} + (\nabla \boldsymbol{u})^{\mathsf T}), & \text{in} \Omega.
		\end{cases}
	\end{aligned}
\end{equation}
This misspecified prior neglects the nonlinear deformation effects in the hyperelastic model and therefore introduces a systematic discrepancy between the prescribed physics and the reference solutions. For data generation, we sample $300$ compression levels with $\epsilon$ ranging from $0$ to $0.2$, and compute the corresponding displacement solutions by a finite element solver. The nonlinear system is solved by Newton's method with line search, using GMRES for the linear subproblems together with a continuation strategy. The dataset is split into $200$ training samples and $100$ test samples. 

During training, the values of the right-boundary displacement function at $21$ equidistant points are used as the input to the solution operator. These boundary measurements are also concatenated with the displacement field predicted by the solution operator on a regular $51\times21$ grid, thereby forming the input to the correction operator. For boundary observations, $101$ equally spaced points are specified on both the top and bottom boundaries, while $21$ equally spaced points are prescribed on each of the left and right boundaries. In the interior of the domain, we assume that $200$ observations of the displacement field are available. The collocation points for enforcing the governing equation are chosen as the $101\times21$ interior grid points, and all final predictions are evaluated on the same grid. For the weighting coefficients associated with the individual loss terms, we set $\lambda_{\text{e}} = 100, \lambda_{\text{u}} = 100000, \lambda_{\text{bc}}=1$.

\cref{fig:2d-hyper-u} presents the prediction results of the displacement field components for a representative test sample with $\epsilon=0.0991$. For both $u_x$ and $u_y$, the misspecified model captures the main spatial trend of the reference solution, but still exhibits visible local discrepancies. By contrast, the corrected model provides displacement predictions that are closer to the reference solution and shows reduced absolute errors in both components. These results indicate that the learned correction term can effectively compensate for the systematic bias caused by the misspecified linear governing equation. The quantitative results are summarized in \cref{tab:2d-hyper-error}, which reports the relative $L_2$ errors on the test set for the four models considered above. Although the visual differences between the misspecified and corrected models may not always be pronounced in the solution plots, the quantitative errors reveal a clear improvement after applying the correction strategy. Although the standard DeepONet improves over the misspecified model by relying on the available data, its accuracy remains inferior to that of the corrected model under the same limited observation setting. As expected, the known model achieves the smallest errors and serves as the reference upper benchmark among the compared methods.

\begin{table}[htbp]
	\centering
	\caption{\textbf{2d hyperelastic problem.} Relative $L_2$ errors on the test set for displacement components $u_x$, $u_y$, and $|\boldsymbol{u}|$, computed with respect to the reference solution for the known model, misspecified model, corrected model, and standard DeepONet. Values are reported as mean $\pm$ standard deviation over five independent runs. All values are in units of $10^{-3}$.}\label{tab:2d-hyper-error}	
	\begin{tabular}{lccc}
		\toprule
		Model & $u_x$ & $u_y$ & $|\boldsymbol{u}|$ \\
		\midrule
		
		Known model
		& $0.68 \pm 0.17$
		& $0.61 \pm 0.02$
		& $0.46 \pm 0.04$ \\
		
		Misspecified model
		& $25.50 \pm 0.09$
		& $11.00 \pm 0.04$
		& $11.51 \pm 0.04$ \\
		
		Standard DeepONet
		& $6.45 \pm 0.34$
		& $5.83 \pm 0.51$
		& $4.19 \pm 0.27$ \\
		
		Corrected model
		& $3.75 \pm 0.20$
		& $0.97 \pm 0.26$
		& $2.06 \pm 0.13$ \\
		
		\bottomrule
	\end{tabular}
\end{table}

\section{Conclusion}\label{sec:conclusion}
We have proposed a unified operator correction framework for learning solution operators in the presence of misspecified governing PDEs. The proposed framework formulates the correction at the operator level, where an approximate physical model provides a prior operator and a learnable correction operator is introduced to compensate for the discrepancy caused by model misspecification. In this work, the framework is implemented using a serial DeepONet architecture. The first DeepONet produces a baseline prediction associated with the solution function, while the second DeepONet learns a correction conditioned on both the input function and the solution prediction. The learned correction is then incorporated into the physics constraint, enabling the solution operator and the correction operator to be trained jointly. Numerical results on several benchmark problems show that the proposed method effectively mitigates the adverse effects of incorrect physical constraints, improves predictive accuracy, and enhances training stability compared with the misspecified model. Experiments with sparse and noisy observations further indicate that the framework remains robust under imperfect data. Comparisons with standard DeepONet also show that incorporating an approximate physical prior together with a learned correction can more effectively exploit limited observational data. Although DeepONet is adopted as the main implementation in this study, the proposed correction principle is architecture independent and can be extended to other neural operator models. Future work will consider applications to time-dependent multiphysics systems, inverse problems, and real experimental data.

%\section*{Acknowledge}
%The first author is partially supported by 

\bibliography{ref.bib}

\end{document}